\documentclass[11pt]{article}
\usepackage{amsfonts}
\def\bC{\mathbf{\overline{C}}}
\def\R{\mathbf{R}}

\def\C{\mathbf{C}}

\def\T{\mathbf{T}}
\def\dist{\mathrm{dist}}

\def\diam{\mathrm{diam}}

\newtheorem{theorem}{Theorem}
\newtheorem{lemma}{Lemma}

\title{On the length of lemniscates}
\date{\,}
\author{Alexandre 
Eremenko\footnote{Supported by EPSRC grant GR/L 35546
at Imperial College and by  NSF grant DMS-9800084}{$\;\;$} 
and Walter Hayman}
\begin{document}
\maketitle
\begin{abstract}
We show that for a polynomial $p(z)=z^d+\ldots$ the length of the
level set $E(p):=\{ z:|p(z)|=1\}$ is at most $9.173\, d$, which improves an
earlier estimate due to P. Borwein. For $d=2$ we show that the
extremal level set is the Bernoullis' Lemniscate. One ingredient of our 
proofs is the fact that for an extremal polynomial the set $E(p)$ is connected.
\end{abstract}

For a monic polynomial $p$ of degree $d$ we write $E(p):=\{ z:|p(z)|=1\}$.
A conjecture of Erd\H{o}s, Herzog and Piranian \cite{EHP}, repeated by
Erd\H{o}s
in \cite{Hayman} and elsewhere, is that the length $|E(p)|$ is maximal
when $p(z):=z^d+1$. It is easy to see that in this conjectured
extremal case $|E(p)|=2d+O(1)$
when $d\to\infty$.

The first upper estimate $|E(p)|\leq 74d^2$ is due to Pommerenke \cite{Pommerenke}.
Recently P. Borwein \cite{Borwein}
gave an estimate which is linear in $d$, namely
$$|E(p)|\leq 8\pi e d\approx 68.32 d.$$ Here we improve Borwein's result.

Let $\alpha_0$ be the least upper bound of perimeters of the convex hulls of
compact connected sets of logarithmic capacity $1$. The precise value of
$\alpha_0$ is not known, but Pommerenke \cite{Ped2} proved the estimate
$\alpha_0<9.173$. The conjectured value is $\alpha_0=3^{3/2}2^{2/3}\approx 8.24$.

\begin{theorem} For monic polynomials $p$ of degree $d$ $|E(p)|\leq \alpha_0 d<
9.173 d$.
\end{theorem}

A similar problem for rational functions turns out to be much easier,
and can be solved completely by means of Lemma \ref{lemma4} below.

\begin{theorem} Let $f$ be a rational function of degree $d$. Then the
spherical
length of the preimage under $f$ of any circle $C$ is at most
$d$ times the length of a great circle.
\end{theorem}
This is best possible as the example $f(z)=z^d$ and $C=\R$ shows.

{\em Remarks.} Borwein notices that his method would give the estimate
$4\pi d\approx 12.57 d$ if one knew one of the following facts:
a) the precise estimate of the size of the exceptional set
in Cartan's Lemma (Lemma \ref{cartan} below) or
b) the fact that for extremal polynomials the set $E(p)$ is
connected. The second fact turns out to be correct (this is our
Lemma~\ref{cartan}), and in addition we
can improve from $4\pi$ to $9.173$ by
using more precise arguments than those
of Borwein.

The main property of the level sets $E(p)$ is the following
\begin{lemma}\label{lemma4} For every rational function $f$ of degree $d$ the
$f$-preimage of any line or circle has no more than $2d$ intersections with any
line or circle $C$, except finitely many $C$'s.
\end{lemma}

{\em Proof}. The group of fractional-linear transformations acts transitively
on the set of all circles on the Riemann sphere, and a composition of
a rational function with a fractional-linear transformation is a rational
function of the same degree.

Thus it is enough to prove that for a rational function $f$ of degree $d$
the set $F:=\{ z:f(z)\in\R\}$ has at most $2d$ points of intersection
with the real line $\R$, unless $\R\subset F$.
Let $z_0$ be such a point of intersection.
Then $z_0$ is a zero of the rational function
$f_1(z)=f(z)-\overline{f(\overline{z})}$. But $f_1$ evidently has degree at
most $2d$ and thus cannot have more than $2d$ zeros, unless $f_1\equiv 0$.

\hfill$\Box$
\medskip

The length of sets described in Lemma \ref{lemma4} can be estimated using
the following lemma, in which we denote
by $\pi_x$ and $\pi_y$ the orthogonal projections 
onto a pair of perpendicular coordinate
axes.

\begin{lemma}\label{lemma5} If an analytic curve $\Gamma$ intersects each vertical and
horizontal line at most $n$ times than
$|\Gamma|\leq n(|\pi_x(\Gamma)|+|\pi_y(\Gamma)|)$. 
\end{lemma}

{\em Proof}.
We break the curve $\Gamma$ into finitely many pieces $l_j$ such that every
$l_j$ intersects each vertical or horizontal line at most once. Then
we have $|l_j|\leq|\pi_x(l_j)|+|\pi_y(l_j)|$. We obtain this by
approximating $l_j$ by broken lines whose segments are parallel to the
coordinate axes. Adding these inequalities for all pieces and using the
fact that both projection maps are at most $n$-to-$1$ on $\Gamma$ we 
obtain the result.
\hfill$\Box$
\medskip

\noindent
{\bf Corollary} {\em Every connected subset $l$ of $E(p)$ has the property}
$$|l|\leq 2d(|\pi_x(l)|+|\pi_y(l)|)\leq 4d\,\diam(l).$$

\begin{lemma}\label{cartan} {\rm (H. Cartan, see for example \cite[p. 19]{L}).}
For a monic polynomial $p$ of degree $d$ the set $\{ z:|p(z)|<M\}$ is contained
in the union of discs the sum of whose radii is $2eM^{1/d}$.
\end{lemma}

Pommerenke \cite[Satz 3]{Ped} improved the constant $2e$ in this lemma to
$2.59$ but we will not use this result. 


Now we can prove the existence of extremal polynomials for our problem.

\begin{lemma}\label{existence} The length $|E(p)|$ is a continuous
function of the coefficients of $p$. For every positive integer $d$ there
exists a monic polynomial $p_d$ with the property
$|E(p_d)|\geq |E(p)|$ for every
monic polynomial $p$ of degree $d$.
\end{lemma}

{\em Proof}. Every monic polynomial of degree $d$ can be
written as
$$p(z)=\prod_{j=1}^d(z-z_j).$$
We consider vectors $Z=(z_1,\ldots,z_d)$ in $\C^d$ and denote by $p_Z$
the monic polynomial with the zero-set $Z$.

First we show that $|E(p)|\to 0$ as $\diam Z\to \infty,\; p=p_Z$.
Let $M$ be a number such that $M>(4e)^d$. If the diameter of the set
$Z$ is greater than $4Md$ then we can split $Z$ into
two parts, $Z=Z_1\cup Z_2$ such that $\dist(Z_1,Z_2)>4M$.

Indeed, Let $D$ be the union of closed discs of radii $2M$ centered at the
points $z_1,\ldots,z_d$. If $D$ is connected,
then $\diam D\leq 4Md$, contradicting
our assumption. Thus $D$ is disconnected that is $D=D_1\cup D_2$,
where $D_1$ and $D_2$ are disjoint compact sets, and we set $Z_i=Z\cap D_i$
for $i=1,2$, which proves our assertion. \hfill$\Box$

Consider
two polynomials
$$p_k(z):=\prod_{w\in Z_k}(z-w),\quad k=1,2;\quad\mbox{so that}\quad p=p_1p_2.$$
By Lemma \ref{cartan} the union $L$ of two sets
$L_k:=\{ z:|p_k(z)|<M^{-1}\},\;k=1,2,$
can be covered by discs the sum of whose radii is
$4eM^{-1/d}<1$. Thus the sum of the lengths of the projections of $L$ 
satisfies 
\begin{equation}
\label{slp}
|\pi_x(L)|+|\pi_y(L)|\leq 16eM^{-1/d}.
\end{equation}
On the other hand,
each component of a set $L_k$ contains a zero of $p_k$ and has diameter
less than 2 so that 
$\dist(L_1,L_2)>4M-4>2M,$ since $M>4e$.

Next we show that $E(p)\subset L_1\cup L_2$.
Indeed, suppose that $z\in E(p)$. Assume without loss of generality that 
$\dist(z,L_1)\leq\dist(z,L_2).$ Then\newline $\dist(z,L_2)>M$
and thus $|p_2(z)|>M$, so that $$|p_1(z)|=|p(z)|/|p_2(z)|<M^{-1}$$ and this implies
that $z\in L_1$.

We conclude that $|\pi_x(E(p))|+|\pi_y(E(p))|\leq 16eM^{-1/d}$, which tends to $0$
as $M\to\infty$. Now an application of the Corollary after Lemma \ref{lemma5}
concludes the proof of our assertion that $|E(p)|\to 0$ as
$\diam Z\to\infty$.

Now we show that $|E(p_Z)|$ is a continuous function of the vector
\newline $Z=(z_1,\ldots,z_d)\in\C^d$.
Consider the multivalued algebraic function
$$q(Z,w)=(d/dw)(p^{-1}(w)).$$ The coefficients  of the algebraic equation
defining this function $q$ are polynomials of $Z$ and $w$, and $q(Z,w)\neq 0$
in $\C^d\times\C$ because this is a derivative of an inverse function. So
all branches
of $q$ are continuous with respect to $w$ and $Z$ at every
point where these branches are finite
(see, for example \cite[Theorem 12.2.1]{Hille}).
Denoting by $\T$ the unit circle,
we have
$$|E(p)|=\int_{\T}Q(Z,w)\,|dw|,\quad\mbox{where}\quad Q(Z,w)=\sum |q(Z,w)|,$$
and the summation is over all values of the multi-valued function $q$.
To show that this integral is a continuous
function of the parameter $Z$, we will verify that the family of functions
$w\mapsto\sum |q(Z,w)|,\;\T\to\R$ has a uniform integrability property. 

Let $K$ be an arc of the unit circle of length $\delta<\pi/6$.
Then this arc is contained in a disc $D(w,r)$ of radius $r=\delta/2$,
centered at the middle point $w$ of the arc $K$.
By  Lemma \ref{cartan}, applied to $p-w$,
the full preimage $p^{-1}D(w,r)$ can be covered by discs the sum of whose
radii
is at most $2er^{1/d}$. So the sum of the vertical and horizontal projections of
$p^{-1}D(w,r)$ is at most $8er^{1/d}$. Finally by the Corollary after 
Lemma~\ref{lemma5}, the length of the part of $E(p)$ which is mapped to $K$ is
at most $\epsilon:=16der^{1/d}=16de(\delta/2)^{1/d}$. Thus 
\begin{equation}
\label{epsilon}
\int_KQ(Z,w)\,|dw|<\epsilon,
\end{equation}
where $\epsilon\to 0$ as $\delta\to 0$ uniformly with respect to
$Z$.

Suppose now that $Z_0\in\C^d$. Consider the points $w_1,\ldots,w_k$
on the unit circle $\T$, such that $Q(Z_0,w_j)=\infty$. Then $k\leq d-1$,
because a polynomial $p$ of degree $d$ can have no more than $d-1$
critical points.
Given that $\epsilon>0$ we choose open arcs $K_j$ so that $w_j\in K_j\subset\T,\;
1\leq j\leq k$, and (\ref{epsilon}) is satisfied with $K=\cup_jK_j$ whenever
$Z\in\C^d$. Now we have $Q(Z,w)\to Q(Z_0,w)$ as $Z\to Z_0$
 uniformly with respect to $w$ in $\T\backslash K$, so that
$$\left|\int_\T Q(Z,w)\,|dw|-\int_\T Q(Z_0,w)\,|dw|\right|\leq 3\epsilon,$$
when $Z$ in $\C^d$ is close enough to $Z_0$.

We have proved that $Z\mapsto|E(p_Z)|$ is a continuous function in $\C^d$
and that
$|E(p_Z)|\to 0$ as $Z\to\infty$. It follows that a maximum of $|E(p)|$
exists. 

To show that $|E(p)|$ is a continuous function of the coefficients we again
refer to the well-known fact \cite[Theorem 12.2.1]{Hille} that the zeros of a monic polynomial
are continuous functions of its coefficients.
\hfill$\Box$
\medskip

In what follows we will call extremal any polynomial
$p$ which maximises $|E(p)|$ in the set of all monic polynomials of degree $d$.

\begin{lemma}\label{lemma1} There exists an extremal polynomial $p$, such that
all critical points of $p$ 
are contained in $E(p)$.
\end{lemma}

{\em Remarks}. From this lemma it follows that the polynomial $z^2+1$ is
extremal for $d=2$.  The level set $\{ z:|z^2+1|=1\}$ is known as
Bernoullis' Lemniscate (it is also one of Cassini's ovals) 
and its length is expressed by the elliptic
integral
$$2^{3/2}\int_0^1\frac{1}{\sqrt{1-x^4}}dx\approx 7.416\,.$$ 
This curve as well as the integral played an important role in the history
of mathematics,  see for example \cite{Pra}.

{\em Proof of Lemma \ref{lemma1}}.

Let $p$ be a polynomial and $a$ a critical value of $p$,
such that 
$a$ does not lie on the unit circle $\T$. Let $U$ be an 
open disc centered at $a$, such that $U$ 
does not contain other critical values. Let
$\Phi:\C\to\C$
be a smooth function whose support is contained in $U$ and such that
$\Phi(a)=1$.
If $\lambda\in\C$ and $\lambda$ satisfies $|\lambda|<\epsilon:=(\max_U|{\mbox{grad}}\Phi|)^{-1}$
then the map
$\phi_\lambda:\C\to\C,\;
\phi_\lambda(z)=z+\lambda\Phi(z)$, is a smooth quasiconformal 
homeomorphism of $\C$. So we have a family of quasiconformal 
homeomorphisms, depending analytically on $\lambda$ for $|\lambda|<\epsilon$.

The composition $q_\lambda:=\phi_\lambda\circ p$ is a family of
quasiregular maps of the plane into itself. We denote by $\mu_f$ the
Beltrami coefficient of a quasiregular map $f$, that is
$\mu_f:=f_{\overline{z}}/f_z$, where $f_z:=\partial f/\partial z$ and
$f_{\overline{z}}:=\partial f/\partial\overline{z}$.
By the chain rule (see for example \cite[p. 9]{Ahlfors})
\begin{equation}\label{chain}
\mu_{q_\lambda}=
\left(\frac{|p^\prime|}{p^\prime}\right)^2 \mu_{\phi_\lambda}\circ p,
\end{equation}
so that $\mu_{q_\lambda}$
depends analytically on $\lambda$ for $|\lambda|<\epsilon$. 
According to the Existence and Analytic Dependence
on Parameter Theorems for the Beltrami equation
(see, for example \cite[Ch. I, theorems 7.4 and 7.6]{Carleson}),
there exists a family of quasiconformal homeomorphisms $\psi_\lambda:\C\to\C$,
satisfying the Beltrami equations 
$$\mu_{\psi_\lambda}=\mu_{q_\lambda},$$
normalized by $\psi_\lambda=z+o(1),\;z\to\infty$, and analytically depending
on $\lambda$ for $|\lambda|<\epsilon$ for every fixed $z$.

It follows that 
$$p_\lambda:=q_\lambda\circ\psi_\lambda^{-1}=\phi_\lambda\circ p\circ
\psi_\lambda^{-1}$$
are entire functions. As they are all $d$-to-$1$, they are polynomials
of degree $d$, and the normalization of $\psi$ implies that these polynomials
are monic. These polynomials $p_\lambda$ may be considered as obtained from
$p$ by shifting one critical value from $a$ to $a+\lambda$, while all other
critical values remain unchanged. The functions $\lambda\mapsto p_\lambda(z)$
are continuous (in fact analytic) for every $z$. Thus the coefficients
of $p_\lambda$ are continuous functions of $\lambda$. It follows by Lemma
\ref{existence} that  $|E(p_\lambda)|$ is a continuous function of $\lambda$.

Now we assume that $p$ is an extremal polynomial, and that a critical
value $a$ of $p$ does not belong to the unit circle $\T$. Then we can
choose the disc $U$ in the construction above such that $U$ does
not intersect the unit circle. As $\phi_\lambda$ is conformal outside
$U$, we conclude from (\ref{chain}) that $q_\lambda$ and thus $\psi_\lambda$
are conformal away from $p^{-1}(U)$. This implies that 
$\psi_\lambda$ is conformal in the neighborhood of $E(p)$,
and we have
$$|E(p_\lambda)|=|\psi_\lambda(E(p))|=\int_{E(p)}\left|
\frac{d\psi_\lambda}{dz}\right|\,|dz|.$$
As $\phi_\lambda$ depends analytically on $\lambda$ so does $d\psi_\lambda/dz$;
thus $|d\psi_\lambda/dz|$ is a subharmonic function of $\lambda$ for
$|\lambda|<\epsilon$ for every
fixed $z$. It follows that $|E(p_\lambda)|$ is subharmonic
for $|\lambda|<\epsilon$.
Because we assumed that $p$ is extremal, this subharmonic function has a
maximum at the point $0$, so it is constant. 

Now we consider all critical values $a_1,\ldots,a_n$ of $p$ which do not
belong to the unit circle $\T$, and connect each $a_j$ with $\T$ 
by a curve $\gamma_j$ such that all these curves are disjoint and do
not intersect $\T$, except at one endpoint. Performing the deformation
described above, we move all critical values $a_j$, one at a time,
along $\gamma_j$ to
the unit circle, and obtain as a result a monic polynomial $p^*$ of degree $d$,
all of whose
critical values belong to $\T$.   
This is equivalent to the property that all critical points of $p^*$ belong to 
$E(p^*)$. We have $|E(p^*)|=|E(p)|$, because $|E(p)|$ remains constant as
a critical value $a_j$ is moved along $\gamma_j$. Thus $p^*$ is also extremal.
\hfill$\Box$

\begin{lemma}\label{lemma2} There exists an extremal polynomial $p$ for which
the set
$E(p)$ is connected.
\end{lemma}

{\em Proof}. Put $D=\{ z\in\bC:|P(z)|>1\},$ and $\Delta=\{ z\in\bC:|z|>1\}$.
Let $p$ be an extremal polynomial constructed as in Lemma \ref{lemma1}.
Then $p:D\to\Delta$ is a ramified covering of degree $d$ having exactly
one critical
point of index $d-1$, namely the point $\infty$.
By the Riemann--Hurwitz Formula  $D$ is simply connected,
so $E(p)$ is connected. \hfill$\Box$
\medskip

{\em Remarks}. By moving those critical values whose moduli are
greater than $1$ towards
infinity, rather than to the unit circle, and using the arguments from 
the proof of Lemma 4, one can show that an extremal
polynomial cannot have critical
values with absolute value greater than $1$. It follows that in fact for all
extremal polynomials $p$ the level sets $E(p)$ are connected.
We will not use this
additional information in the proof of Theorem 1. 

\begin{lemma}\label{lemma3} {\rm (Pommerenke \cite[Satz 5]{Ped2})}.
Let $E$ be a connected compact set of logarithmic capacity $1$.
Then the perimeter of the convex hull of $E$ is at most
$$\pi(\sqrt{10}-3\sqrt{2}+4)<9.173.$$
\end{lemma}
\medskip

{\em Proof of Theorem 1}. Let $p$ be an extremal polynomial with connected
set $E(p)$. Such a $p$ exists by Lemma \ref{lemma2}.
Applying Lemma \ref{lemma3}
we conclude that the perimeter of the convex hull of $E$ is at most
$9.173$. 

Now the integral-geometric formula \cite{Santalo} for the length of a curve gives
$$|E|=\frac{1}{2}\int_0^\pi\int_{-\infty}^{\infty}N_E(\theta,x)\,dx\,d\theta,$$
where $N_E(\theta,x)$ is the number of intersections of $E$ with the
line $$\{ z:\Re(ze^{-i\theta})=x\}.$$ A connected compact set $E$ intersects exactly
those lines which the boundary of its convex hull intersects.
But the boundary of the convex hull intersects almost every line
either $0$ or $2$ times, while a set $E(p)$ intersects
each line at most $2\deg p$ times. Thus $|E|<9.173 d$. This proves our assertion.
\medskip

{\em Proof of Theorem 2}. Following Borwein we use the Poincar\'{e}
Integral-Geometric
Formula \cite{Santalo}. Assuming that the great circles have length $2\pi$, we denote by
$l(E)$ the spherical length of $E$, by $dx$ the spherical area
element and by $v(E,x)$ the number of intersections
of $E$ with the great circle, one of whose centers is $x$.
The  Poincar\'{e} Formula is
$$l(E)=\frac{1}{4}\int v(E,x)\,dx.$$
Now if $E(f)$ is the preimage of a circle under a rational
function $f$ of degree $d$ then by Lemma \ref{lemma4} $E(f)$ intersects every
great circle at most $2d$ times, so that the spherical length $l(E(f))$ is
at most $2\pi d$.\hfill$\Box$
\medskip

We are very grateful to Christian
Pommerenke for helpful discussion and references. We also thank the referee,
whose suggestions improved our original estimate in Theorem 1.

\medskip

\noindent
Purdue University, West Lafayette IN 47907
\newline
{\em eremenko@math.purdue.edu}
\medskip

\noindent
Imperial College, London SW7 2BZ
\medskip

\end{document}